\documentclass[psamsfonts,reqno]{amsart}
\usepackage{amssymb,eucal,latexsym,xy}
\xyoption{all}
\xyoption{ps}

\theoremstyle{plain}

\newtheorem{lemma}{Lemma}[section]
\newtheorem{theorem}[lemma]{Theorem}

\newtheorem{corollary}[lemma]{Corollary}
\newtheorem{claim}{Claim}

\newtheorem*{stat}{\name}
\newcommand{\name}{testing}

\theoremstyle{definition}
\newtheorem{definition}[lemma]{Definition}

\newtheorem*{problem}{Problem}

\theoremstyle{remark}

\newtheorem{notation}[lemma]{Notation}

\newenvironment{all}[1]{\renewcommand{\name}{#1}\begin{stat}}
                         {\end{stat}}

\newcommand{\qedc}{{\qed}~{\rm Claim~{\theclaim}.}}

\newenvironment{cproof} {\begin{proof}[Proof of Claim.]}
{\qedc\renewcommand{\qed}{}\end{proof}}

\numberwithin{equation}{section}

\newcommand{\pup}[1]{\textup{(}{#1}\textup{)}}

\newcommand{\ba}{{\boldsymbol{a}}}
\newcommand{\bb}{{\boldsymbol{b}}}
\newcommand{\bc}{{\boldsymbol{c}}}
\newcommand{\bd}{{\boldsymbol{d}}}
\newcommand{\be}{{\boldsymbol{e}}}
\newcommand{\xe}{{\mathbf{e}}}

\newcommand{\bx}{{\boldsymbol{x}}}
\newcommand{\by}{{\boldsymbol{y}}}
\newcommand{\bz}{{\boldsymbol{z}}}

\newcommand{\Pow}{\mathfrak{P}}

\newcommand{\Vhom}{V-ho\-mo\-mor\-phism}

\newcommand{\cm}{com\-mu\-ta\-tive mo\-no\-id}
\newcommand{\poag}{partially ordered abelian group}
\newcommand{\ppoag}{pointed \poag}

\newcommand{\es}{\varnothing}
\newcommand{\into}{\hookrightarrow}
\newcommand{\onto}{\twoheadrightarrow}
\newcommand{\ip}{\mathbin{\bowtie}}

\newcommand{\set}[1]{\{#1\}}
\newcommand{\setm}[2]{\set{#1\mid#2}}
\newcommand{\seq}[1]{({#1})}
\newcommand{\seqm}[2]{\seq{#1\mid#2}}
\newcommand{\famm}[2]{\left(#1\mid#2\right)}

\newcommand{\ol}[1]{\overline{#1}}

\newcommand{\go}{\omega}

\DeclareMathOperator{\supp}{supp}
\DeclareMathOperator{\FP}{FP}

\newcommand{\cD}{\mathcal{D}}
\newcommand{\cL}{\mathcal{L}}

\newcommand{\ZZ}{\mathbb{Z}}
\newcommand{\FF}{\mathbb{F}}

\DeclareMathOperator{\im}{im}
\DeclareMathOperator{\rng}{rng}

\DeclareMathOperator{\Conc}{Con_c}

\DeclareMathOperator{\Dim}{Dim}

\DeclareMathOperator{\DD}{\Delta}
\newcommand{\xI}{\mathbf{I}}
\newcommand{\xE}{\mathbf{E}}
\newcommand{\xF}{\mathbf{F}}

\newcommand{\id}{\mathrm{id}}
\newcommand{\jz}{$(\vee,0)$}
\newcommand{\jzs}{\jz-semi\-lat\-tice}

\newcommand{\fin}[1]{[#1]^{<\go}}

\newcommand{\fact}[2]{{#1}_{[{#2}]}}

\begin{document}

\title[A dimension group with order-unit]%
{A $K_0$-avoiding dimension group with\\
an order-unit of index two}
\author[F.~Wehrung]{Friedrich~Wehrung}
\address{LMNO, CNRS UMR 6139\\
         D\'epartement de Math\'ematiques, BP 5186\\
         Universit\'e de Caen, Campus 2\\
         14032 Caen Cedex\\
         France}
 \email{wehrung@math.unicaen.fr}
 \urladdr{http://www.math.unicaen.fr/\~{}wehrung}

 \date{\today}

\begin{abstract}
We prove that there exists a dimension group $G$ whose positive cone is not
isomorphic to the dimension monoid $\Dim L$ of any lattice $L$. The
dimension group~$G$ has an order-unit, and can be taken of any cardinality
greater than or equal to $\aleph_2$. As to determining the
positive cones of dimension groups in the range of the $\Dim$ functor, the
$\aleph_2$ bound is optimal. This solves negatively the problem, raised by
the author in 1998, whether any conical refinement monoid is isomorphic to
the dimension monoid of some lattice. Since $G$ has an order-unit of
index~$2$, this also solves negatively a problem raised in 1994 by K.\,R.
Goodearl about representability, with respect to $K_0$, of dimension groups
with order-unit of index~$2$ by unit-regular rings.
\end{abstract}

\maketitle

\section*{Introduction}\label{S:Intro}
The nonstable $K$-theory of a ring $R$ studies the category of finitely
generated projective right $R$-modules. The lattice-theoretical analogue
of nonstable K-theory is encoded by the \emph{dimension monoid} functor.
The dimension monoid of a lattice~$L$ (see \cite{WDim}) is the \cm\
defined by generators $\DD(x,y)$, for $x\leq y$ in~$L$, and relations
 \begin{itemize}
 \item[(D0)] $\DD(x,x)=0$, for all $x\in L$.
 \item[(D1)] $\DD(x,z)=\DD(x,y)+\DD(y,z)$, for all $x\leq y\leq z$ in $L$.
 \item[(D2)] $\DD(x\wedge y,x)=\DD(y,x\vee y)$, for all $x$, $y\in L$.
 \end{itemize}
The dimension monoid $\Dim L$ is a precursor of the
semilattice $\Conc L$ of compact congruences of $L$, in the sense that
$\Conc L$ is isomorphic to the maximal semilattice quotient of $\Dim L$,
see \cite[Corollary~2.3]{WDim}. Furthermore, although it is still an open
problem whether $\Dim L$ is a refinement monoid (see
Section~\ref{S:Basic} for a definition) for every lattice $L$ (see
\cite[Problem~3]{WDim}), the answer is known for a few large classes of
lattices, namely, the class of all modular lattices
(\cite[Theorem~5.4]{WDim}) and the class of all lattices without infinite
bounded chains (see Theorem~6.18 and Corollary~7.8 in \cite{WDim}).

The question of a converse, namely whether every refinement monoid is
isomorphic to the dimension monoid of some lattice, was raised by the
author in \cite[Problem~4]{WDim}. This question is an analogue, for the
$\Dim$ functor, of the Congruence Lattice Problem that asks whether every
distributive \jzs\ is isomorphic to $\Conc L$, for some lattice $L$ (see
\cite{CLPSurv} for a survey). Partial positive answers were known. For
example, it follows from \cite[Theorem~1.5]{GoHa86} and results
in~\cite{WDim} (see the proof of Corollary~\ref{C:MainVNeg}) that for
every dimension group $G$ of cardinality at most~$\aleph_1$, the positive
cone~$G^+$ is isomorphic to~$\Dim L$ for some sectionally complemented,
modular lattice~$L$. For the cardinality~$\aleph_2$ and above, the problem
was still open. Different, though related, positive results about the
dimension theory of complete modular lattices but also of self-injective
modules or AW*-algebras, are established in~\cite{Gore}. In
particular, the dimension monoids of complete, complemented,
modular, upper continuous lattices are completely characterized.

\begin{all}{Main theorem}
There exists a dimension group $G$ with order-unit of index~$2$ such that
for any lattice $L$, the positive cone $G^+$ of $G$ is not the image of
$\Dim L$ under any \Vhom. Furthermore, $G$ may be taken of any cardinality
greater than or equal to $\aleph_2$.
\end{all}

(We refer to Section~\ref{S:Basic} for precise definitions.)
In particular, $G^+$ is not isomorphic to $\Dim L$, for any lattice $L$.
This solves \cite[Problem~4]{WDim}. Also, \emph{$G$ is not isomorphic to
$K_0(R)$, for any unit-regular ring $R$} (see Corollary~\ref{C:MainVNeg}),
which solves negatively the problem raised by K.\,R. Goodearl on the last
page of \cite{Good94}. A stronger and more precise statement of the main
theorem is presented in Theorem~\ref{T:MainVNeg}.

The proof of our result is based on the proofs of earlier
counterexamples, the first of this sort, due to the author in
\cite{NonMeas}, being a dimension group with order-unit of cardinality
$\aleph_2$ that is not isomorphic to $K_0(R)$, for any von~Neumann regular
ring $R$. Later counterexamples to related questions in lattice theory
appeared in \cite{UnifRef,PTW,TuWe}. A common point of their proofs is
that they all use the \emph{Kuratowski Free Set Theorem}, in the form of
Lemma~\ref{L:Kura}. Also, they all express that certain distributive
semilattices cannot be expressed as $\Conc L$, for lattices $L$ with
\emph{permutable congruences}.

By contrast, the proof of our main theorem does not require any assumption
about permutable congruences on the lattice $L$. Also, unlike the
construct of~\cite{NonMeas}, our dimension group counterexample is not a
rational vector space. This is also the case for the dimension groups
considered in~\cite{Good94}, in which the order-unit has \emph{finite
index}. However, in \cite{Good94} are proven \emph{positive} results, not
from the viewpoint of the dimension theory of lattices but from the
closely related (see Lemma~\ref{L:V(R)DimL}) viewpoint of the
nonstable K-theory of von~Neumann regular rings. For
example \cite[Theorem~4.3]{Good94}, whenever $G$ is an abelian
lattice-ordered group with order-unit of finite index, there exists a
biregular locally matricial algebra $R$ such that
$G\cong\nobreak K_0(R)$; hence $G^+\cong\Dim L$, where $L$ is the lattice
of all principal right ideals of~$R$, see the proof of
Corollary~\ref{C:MainVNeg} (as $R$ is unit-regular, it is sufficient to
use~$R$ instead of~$M_2(R)$).

\section{Basic concepts}\label{S:Basic}

Every \cm\ will be endowed with its \emph{algebraic quasi-ordering},
defined by
 \[
 x\leq y\ \Longleftrightarrow\ (\exists z)(x+z=y).
 \]
We say that $M$ is \emph{conical}, if $x\leq 0$ implies that $x=0$, for
all $x\in M$. For \cm s $M$ and $N$, a monoid homomorphism
$\mu\colon M\to N$ is a \emph{\Vhom}, if whenever $c\in M$ and $\ol{a}$,
$\ol{b}\in N$ such that $\mu(c)=\ol{a}+\ol{b}$, there are $a$, $b\in M$
such that $c=a+b$, $\mu(a)=\ol{a}$, and $\mu(b)=\ol{b}$. An \emph{o-ideal}
of a \cm\ $M$ is a nonempty subset $I$ of $M$ such that $x+y\in I$ if{f}
$x$, $y\in I$, for all $x$, $y\in M$. For an o-ideal $I$ of a \cm\ $M$,
the least monoid congruence $\equiv_I$ that identifies all elements of $I$
to zero is defined by
 \[
 x\equiv_Iy\ \Longleftrightarrow\ (\exists u,v\in I)(x+u=y+v),
 \quad\text{for all }x,\,y\in M.
 \]
We denote by $M/I$ the quotient monoid $M/{\equiv_I}$, and we denote by
$[x]_I$ the $\equiv_I$-class of any $x\in M$. The proof of the following
lemma is straightforward.

\begin{lemma}\label{L:Vhomzs}
Let $M$ and $N$ be \cm s with $N$ conical and let $\mu\colon M\to N$ be a
monoid homomorphism. Then the subset $I=\setm{x\in M}{\mu(x)=0}$ is an
o-ideal of $M$, and there exists a unique monoid homomorphism
$\ol{\mu}\colon M/I\to N$ such that $\ol{\mu}([x]_I)=\mu(x)$ for all
$x\in M$. Furthermore, if $\mu$ is a \Vhom, then so
is~$\ol{\mu}$.
\end{lemma}

A \cm\ $M$ is a \emph{refinement monoid}, if $a_0+a_1=b_0+b_1$ in~$M$
implies the existence of $c_{i,j}\in M$, for $i$, $j<2$, such that
$a_i=c_{i,0}+c_{i,1}$ and $b_i=c_{0,i}+c_{1,i}$, for all $i<2$. A \jzs\
$S$ is \emph{distributive}, if it is a refinement monoid. Equivalently,
the ideal lattice of $S$ is distributive, see \cite[Section~II.5]{GLT2}.

We use the notation, terminology, and results of \cite{Gpoag} for \poag s.
For \poag s $G$ and $H$, a group homomorphism $f\colon G\to H$ is a
\emph{positive homomorphism}, if $f[G^+]\subseteq H^+$. For a \poag\ $G$
and a positive integer~$n$, we say that an element $e\in G^+$ has
\emph{index at most~$n$}, if
$(n+1)x\leq e$ implies that $x=0$, for all $x\in G^+$. We say that
$e\in G^+$ is an \emph{order-unit} of $G$, if for all $x\in G$, there
exists a natural number $n$ such that $x\leq ne$.

We say that a \poag\ $G$ is
\begin{itemize}
\item[---] an \emph{interpolation group}, if for all $x$, $x'$, $y$,
$y'\in G$, if $x,x'\leq y,y'$, then there exists $z\in G$ such that
$x,x'\leq z\leq y,y'$;

\item[---] \emph{unperforated}, if $mx\geq0$ implies that $x\geq0$, for
every $x\in G$ and every positive integer~$m$;

\item[---] \emph{directed}, if $G=G^++(-G^+)$;

\item[---] a \emph{dimension group}, if $G$ is a directed, unperforated
interpolation group.
\end{itemize}

Particular cases of dimension groups are the \emph{simplicial groups},
that is, the \poag s isomorphic to finite powers of the additive
group~$\ZZ$ of all integers, ordered componentwise. A theorem of Effros,
Handelman, and Shen states that dimension groups are exactly the direct
limits of simplicial groups, but we shall not need this result in the
present paper.

A \emph{\ppoag} is a pair $\seq{G,e_G}$, where $G$ is a \poag\ and
$e_G\in G^+$. We shall call $e_G$ the \emph{distinguished element} of
$\seq{G,e_G}$. For \ppoag s $\seq{G,e_G}$ and $\seq{H,e_H}$, a positive
homomorphism $f\colon G\to H$ is \emph{normalized}, if $f(e_G)=e_H$. We
shall write \ppoag s either in the form $\seq{G,e_G}$ in case the
distinguished element $e_G$ needs to be specified, or simply $G$ otherwise.

For any lattice $L$, the symbol $\DD({}_-,{}_-)$ is extended to any pair
of elements of $L$, by defining $\DD(x,y)=\DD(x\wedge y,x\vee y)$, for
all $x$, $y\in L$. The map $\DD$ thus extended satisfies all the basic
properties defining distances, see \cite[Proposition~1.9]{WDim}.

\begin{lemma}\label{L:DDDist}
The following statements hold, for all $x$, $y$, $z\in L$:
\begin{enumerate}
\item $\DD(x,y)=0$ if{f} $x=y$;
\item $\DD(x,y)=\DD(y,x)$;
\item $\DD(x,z)\leq\DD(x,y)+\DD(y,z)$.
\end{enumerate}
\end{lemma}

Of course, in (iii) above, the \cm\ $\Dim L$ is endowed with its algebraic
quasi-ordering.

The following result is an immediate consequence of
\cite[Lemma~4.11]{WDim}, applied to the partial semigroup of closed
intervals of $L$ endowed with projectivity as in \cite[Section~5]{WDim}.
It concentrates most of the nontrivial information that we will need about
the dimension monoid.

\begin{lemma}\label{L:DecDD(a,b)}
Let $L$ be a modular lattice, let $u\leq v$ in $L$, and let $\ba$,
$\bb\in\Dim L$. If $\ba+\bb=\DD(u,v)$, then there are a positive integer
$n$ and a decomposition $u=w_0\leq w_1\leq\cdots\leq w_{2n}=v$ such that
 \[
 \ba=\sum\famm{\DD(w_{2i},w_{2i+1})}{i<n}\text{ and }
 \bb=\sum\famm{\DD(w_{2i+1},w_{2i+2})}{i<n}.
 \]
\end{lemma}

For a unital ring $R$, we denote by $\FP(R)$ the category of all finitely
generated projective right $R$-modules, and by $V(R)$ the monoid of all
isomorphism classes of members of $\FP(R)$, see \cite{Good95}.
The monoid $V(R)$ encodes the so-called \emph{nonstable K-theory} of $R$.
If $[X]$ denotes the isomorphism class of a member $X$ of $\FP(R)$, then
the addition of $V(R)$ is defined by $[X]+[Y]=[X\oplus Y]$, for all $X$,
$Y\in\FP(R)$. The monoid $V(R)$ is, of course, always conical. In case $R$
is von~Neumann regular (that is, for all $x\in R$ there exists $y\in R$
such that $xyx=x$), $V(R)$ is a refinement monoid, see
\cite[Theorem~2.8]{GvnRR}.

It is well-known that for a von~Neumann regular ring $R$, the matrix ring
$M_2(R)$ is von~Neumann regular \cite[Theorem~1.7]{GvnRR}. Denote by
$\mathcal{L}(R)$ the (complemented, modular) lattice of principal right
ideals of $R$. The nonstable K-theory of von~Neumann regular rings and the
dimension theory of lattices are related by the following result, which is
an immediate consequence of \cite[Proposition~10.31]{WDim}.

\begin{lemma}\label{L:V(R)DimL}
Let $R$ be a von~Neumann regular ring, and put $L=\mathcal{L}(M_2(R))$.
Then $V(R)\cong\Dim L$.
\end{lemma}

An example due to G.\,M. Bergman, see \cite[Example~4.26]{GvnRR}, shows
that $\mathcal{L}(M_2(R))$ cannot be replaced by $\mathcal{L}(R)$ in the
statement of Lemma~\ref{L:V(R)DimL}.

For a set $X$ and a natural number $n$, we denote by $[X]^n$ (resp.
$[X]^{\leq n}$) the set of all subsets $Y$ of $X$ such that $|Y|=n$
(resp., $|Y|\leq n$). Furthermore, we denote by $\fin{X}$ the
set of all finite subsets of $X$. The set-theoretical core of the proof
of the main theorem consists of the following two results.

\begin{lemma}\label{L:Laza}
Let $X$ be a set of cardinality at least $\aleph_2$ and let
$\Phi\colon X\to\fin{X}$. Then there exists a subset $Y$ of $X$ of
cardinality $\aleph_2$ such that $\eta\notin\Phi(\xi)$, for all distinct
$\xi$, $\eta\in Y$.
\end{lemma}

\begin{proof}
This is a particular case of a result proved by D. L\'az\'ar \cite{Laza}.
See also \cite[Corollary~44.2]{EHMR}.
\end{proof}

\begin{lemma}\label{L:Kura}
Let $X$ be a set of cardinality at least $\aleph_2$, let
$\Psi\colon[X]^2\to\fin{X}$. Then there are distinct $\alpha$, $\beta$,
$\gamma\in X$ such that $\alpha\notin\Psi(\set{\beta,\gamma})$,
$\beta\notin\Psi(\set{\alpha,\gamma})$, and
$\gamma\notin\Psi(\set{\alpha,\beta})$.
\end{lemma}

\begin{proof}
This is a particular case of a result proved by C. Kuratowski
\cite{Kura51}. See also \cite[Theorem~46.1]{EHMR}.
\end{proof}

We denote by $\ZZ^{(X)}$ the additive group of all maps $f\colon X\to\ZZ$
such that the \emph{support} of~$f$, namely $\setm{x\in X}{f(x)\neq0}$, is
finite. A subset $X$ in a partially ordered set $P$ is \emph{cofinal}, if
every element of $P$ lies below some element of $X$.
We identify $n$ with the set $\set{0,1,\dots,n-1}$, for every
natural number $n$.

\section{The functor $\xI$ on \poag s}\label{S:FunctI}

We shall denote by $\cL=\seq{-,0,\xe,\leq,\ip}$ the first-order signature
consisting of one binary operation~$-$ (interpreted as the `difference'),
one binary relation~$\leq$, two constants~$0$ and~$\xe$, and one $4$-ary
operation~$\ip$. Let
$\cD$ denote the class of models of the following axiom system $(\Sigma)$,
written in~$\cL$:
 \[
 (\Sigma):
 \begin{cases}
 \textup{(POAG)}\quad & \text{All axioms of \poag s in }\seq{-,0,\leq}.\\
 \textup{(POINT)}\quad &0\leq\xe.\\
 \textup{(UNPERF)}\quad & \text{Unperforation.}\\
 \textup{(INDEX)}\quad
 &(\forall x)(0\leq 3x\leq\xe\ \Longrightarrow\ x=0).\\
 \textup{(INTERP)}\quad & (\forall x,x',y,y')\bigl(
 x,x'\leq y,y'\ \Longrightarrow\ x,x'\leq\ip(x,x',y,y')\leq y,y'\bigr).\\
 \textup{(SYMM)}\quad & (\forall x,x',y,y')\bigl(
 \ip(x,x',y,y')=\ip(x',x,y,y')=\ip(x,x',y',y)\bigr).
 \end{cases}
 \]
As all axioms of $(\Sigma)$ are universal Horn sentences, it follows from
basic results of the algebraic theory of quasivarieties (see
\cite[Section~V.11]{Malc}) that every model~$G$ for a subsignature~$\cL'$
of~$\cL$ has a unique (up to isomorphism) $\cL'$-homomorphism
$j_G\colon G\to\nobreak\xI(G)$ which is universal among
$\cL'$-homomorphisms from $G$ to some member of $\cD$. This means that
$\xI(G)$ is a member of $\cD$, and for every $\cL'$-homomorphism
$f\colon G\to H$ with $H$ a member of $\cD$, there exists a unique
$\cL$-homomorphism $h\colon\xI(G)\to H$ such that $f=h\circ j_G$.

Applying the universality to the $\cL$-substructure of $\xI(G)$ generated
by the image of $j_G$ yields immediately, in the particular case of
\ppoag s, the following lemma.

\begin{lemma}\label{L:GenI(G)G}
For any \ppoag\ $G$, the structure $\xI(G)$ is the closure, under the
operations $\seq{x,y}\mapsto x-y$ and
$\seq{x,x',y,y'}\mapsto\ip(x,x',y,y')$, of the image of $j_G$.
\end{lemma}

The operation $\ip$ on $\xI(G)$ is a particular instance of the following
notion.

\begin{definition}\label{D:Interpolator}
An \emph{interpolator} on a \poag\ $G$ is a map $\imath\colon G^4\to G$
that satisfies the axioms \textup{(INTERP)} and \textup{(SYMM)} of
the axiom system~$(\Sigma)$. That is,
 \begin{align*}
 & (\forall x,x',y,y'\in G)\bigl(
 x,x'\leq y,y'\ \Longrightarrow\ x,x'\leq\imath(x,x',y,y')\leq
y,y'\bigr).\\
 & (\forall x,x',y,y'\in G)\bigl(
 \imath(x,x',y,y')=\imath(x',x,y,y')=\imath(x,x',y',y)\bigr).
 \end{align*}
\end{definition}

It is obvious that a \poag\ has an interpolator if{f} it is an
interpolation group. We shall naturally view each member of $\cD$ as an
ordered pair $\seq{G,\imath}$, where~$G$ is an unperforated \poag\ and
$\imath$ is an interpolator on~$G$.

For \ppoag s, the meaning of $\xI$ takes the following form: $\xI(G)$ is
a member of $\cD$, the map $j_G$ is a positive homomorphism from~$G$
to~$\xI(G)$, and for every $\seq{H,\imath}\in\cD$ and every
normalized positive homomorphism $f\colon G\to H$, there exists a unique
$\cL$-homomorphism $h\colon\seq{\xI(G),\ip}\to\seq{H,\imath}$ such
that $f=h\circ j_G$. We shall denote this $h$ by $\fact{f}{\imath}$, see
the left hand side diagram of Figure~\ref{Fig:fiota}. In case both~$G$
and~$H$ are \poag s and $f\colon G\to H$ is a normalized positive
homomorphism, the map $\xI(f)=\fact{(j_H\circ f)}{\ip}$ is the unique
$\cL$-homomorphism $h\colon\xI(G)\to\xI(H)$ such that
$h\circ j_G=j_H\circ f$, see the middle diagram of
Figure~\ref{Fig:fiota}.

\begin{figure}[htb]
 \[
 {
 \def\labelstyle{\displaystyle}
 \xymatrix{
 & & & & & & & \seq{G,\imath}\\
 & \seq{H,\imath} & H\ar[r]^{j_H} & \xI(H) & &
 \xI(E)\ar[r]_{\xI(\varphi)}
 \ar[rru]|-(.35){\fact{(f\circ\varphi)}{\imath}} &
 \xI(F)\ar[ru]|-{\fact{f}{\imath}}\\
 G\ar[ru]^f\ar[r]_(.3){j_G} & \seq{\xI(G),\ip}\ar[u]_{\fact{f}{\imath}} &
 G\ar[u]^f\ar[r]_{j_G} & \xI(G)\ar[u]_{\xI(f)} & &
 E\ar[u]^{j_E}\ar[r]_{\varphi} & F\ar[u]^{j_F}\ar[ruu]_f
 }
 }
 \]
\caption{Illustrating $\fact{f}{\imath}$, $\xI(f)$, and
Lemma~\ref{L:Transiota}.}
\label{Fig:fiota}
\end{figure}

Standard categorical arguments give the following two lemmas.
\goodbreak

\begin{lemma}\label{L:DirLimFunc}
The correspondences $G\mapsto\xI(G)$, $f\mapsto\xI(f)$ define a functor
from the category of \ppoag s with normalized positive homomorphisms to the
category $\cD$ with $\cL$-homomorphisms. This functor preserves direct
limits.
\end{lemma}

\begin{lemma}\label{L:Transiota}
Let $E$, $F$, $G$ be \ppoag s, let $\varphi\colon E\to F$ and
$f\colon F\to G$ be normalized positive homomorphisms, and let $\imath$ be
an interpolator on $G$. Then
$\fact{(f\circ\varphi)}{\imath}=\fact{f}{\imath}\circ\xI(\varphi)$
\pup{see the right hand side diagram of Figure~\textup{\ref{Fig:fiota}}}.
\end{lemma}

The following lemma expresses that $\fact{f}{\imath}$ is not `too far'
from $f$.

\begin{lemma}\label{L:Hcliota}
Let $G$ and $H$ be \ppoag s, let $f\colon G\to\nobreak H$ be a normalized
positive homomorphism, and let $\imath$ be an interpolator on $H$.
Then the image of~$\fact{f}{\imath}$ is the least $\imath$-closed subgroup
of $H$ containing the image of~$f$.
\end{lemma}

\begin{proof}
Denote by $H'$ the least $\imath$-closed subgroup
of $H$ containing $\im f$. The subset
$G'=\setm{x\in\xI(G)}{\fact{f}{\imath}(x)\in H'}$ is a subgroup of
$\xI(G)$, closed under the interpolator~$\ip$ as~$H'$ is closed under
$\imath$ and $\fact{f}{\imath}$ is a $\cL$-homomorphism. Since~$G'$
contains~$\im j_G$, it follows from Lemma~\ref{L:GenI(G)G} that
$G'=\xI(G)$.
\end{proof}

The following lemma is even more specific to \ppoag s.

\begin{lemma}\label{L:PresDir}
Let $G$ be a \ppoag. Then the following statements hold:
\begin{enumerate}
\item $\seq{\xI(G),j_G(e_G)}$ is an unperforated pointed interpolation
group with $j_G(e_G)$ of index at most~$2$.

\item The subset $j_G[G]$ is cofinal in $\xI(G)$.

\item If $G$ is directed, then $\xI(G)$ is a dimension group.

\item If $e_G$ is an order-unit of $G$, then $j_G(e_G)$ is an order-unit
of $\xI(G)$.
\end{enumerate}
\end{lemma}

\begin{proof}
(i) is trivial.

(ii) Denote by $H$ the convex subgroup
of $\xI(G)$ generated by the image of $j_G$. Observe that $H$ is closed
under the canonical interpolator $\ip$ of $\xI(G)$,
so it is naturally equipped with a structure of model for $\cL$. Denote by
$f$ the restriction of~$j_G$ from~$G$ to~$H$, and by
$e'\colon H\into\xI(G)$ the inclusion map. Denote by~$h$ the unique
$\cL$-homomorphism from
$\xI(G)$ to $H$ such that $h\circ j_G=f$. {}From $e'\circ h\circ
j_G=e'\circ f=j_G$ and the universal property of $j_G$, it follows that
$e'\circ h=\id_{\xI(G)}$, and so $h(x)=x$, for all $x\in\xI(G)$. Therefore,
$H=\xI(G)$.

(iii) follows immediately from (i) and (ii), while (iv) follows
immediately from~(ii).
\end{proof}

\section{The functors $\xE$ and $\xF$}
\label{S:xExF}

It follows from \cite[Theorem~V.11.2.4]{Malc} that in any quasivariety, one
can form the ``object defined by a given set of generators and relations''.
The following definition uses this general construction in the case of
\ppoag s.

\begin{definition}\label{D:xE(X)}
For a set $X$, we denote by $\seq{\xE(X),\be^X}$ the \ppoag\ defined by
generators $\ba_{\xi}^X$, for $\xi\in X$, and relations
$0\leq\ba_{\xi}^X\leq\be^X$, for $\xi\in X$. We put
$\bb_{\xi}^X=\be^X-\ba_{\xi}^X$, for all $\xi\in X$.
\end{definition}

For $Y\subseteq X$, there are unique positive homomorphisms
$e_{Y,X}\colon\xE(Y)\to\xE(X)$ and $r_{X,Y}\colon\xE(X)\onto\xE(Y)$ such
that
 \begin{align}
 e_{Y,X}(\be^Y)&=\be^X,& e_{Y,X}(\ba^Y_{\eta})&=\ba^X_{\eta},
 \quad\text{for all }\eta\in Y,\label{Eq:DefeYX}\\
 r_{X,Y}(\be^X)&=\be^Y,& r_{X,Y}(\ba^X_{\xi})&=
 \begin{cases}
 \ba^Y_{\xi},&\text{for all }\xi\in Y,\\
 0,&\text{for all }\xi\in X\setminus Y.
 \end{cases}\label{Eq:DefrXY}
 \end{align}
Hence
$r_{X,Y}\circ e_{Y,X}=\id_{\xE(Y)}$, and hence $\xE(Y)$ is a
\emph{retract} of $\xE(X)$. Therefore, we shall identify $\xE(Y)$ with
its image $e_{Y,X}[\xE(Y)]$ in $\xE(X)$, so that $e_{Y,X}$ becomes the
inclusion map from $\xE(Y)$ into $\xE(X)$. Similarly, we shall from now
on write~$\be$ instead of~$\be^X$, $\ba_{\xi}$ instead of $\ba_{\xi}^X$,
and $\bb_{\xi}$ instead of $\bb_{\xi}^X$.

\begin{definition}\label{D:xE(f)}
For sets $X$ and $Y$ and a map $f\colon X\to Y$, we denote by $\xE(f)$
the unique positive homomorphism from $\xE(X)$ to $\xE(Y)$ such that
$\xE(f)(\be)=\be$ and $\xE(f)(\ba_{\xi})=\ba_{f(\xi)}$, for all
$\xi\in X$.
\end{definition}

The proof of the following lemma will introduce a useful explicit
description of the \ppoag\ $\xE(X)$.

\begin{lemma}\label{L:FunctE}
The correspondences $X\mapsto\xE(X)$, $f\mapsto\xE(f)$ define a functor
from the category of sets to the category of all unperforated
\poag s with order-unit. This functor preserves direct limits.
\end{lemma}

\begin{proof}
All items are established by standard categorical arguments, except the
statements about order-unit and, especially, unperforation, that require an
explicit description of $\xE(X)$. Denote by~$\Pow(X)$ the powerset of $X$,
and by $\ol{e}$ the constant function on $\Pow(X)$ with value $1$.
Furthermore, for all $\xi\in X$, we denote by $\ol{a}_{\xi}$ the
characteristic function of $\setm{Y\in\Pow(X)}{\xi\in Y}$. Finally, we let
$F_X$ be the additive subgroup of $\ZZ^{\Pow(X)}$ generated by
$\setm{\ol{a}_{\xi}}{\xi\in X}\cup\set{\ol{e}}$, endowed with its
componentwise ordering. The proof of the following claim is immediate.

\setcounter{claim}{0}
\begin{claim}\label{Cl:DefOrdol}
For all $m\in\ZZ$ and all $\seqm{n_{\xi}}{\xi\in X}\in\ZZ^{(X)}$,
$m\ol{e}+\sum\famm{n_{\xi}\ol{a}_{\xi}}{\xi\in X}\geq0$ in
$F_X$ if{f} $m+\sum\famm{n_{\xi}}{\xi\in Y}\geq0$ in $\ZZ$ for every
$Y\in\Pow(X)$.
\end{claim}

\begin{claim}\label{Cl:GXfreeRel}
There exists an isomorphism from $\xE(X)$ onto $F_X$ that sends $\be$ to
$\ol{e}$ and each $\ba_{\xi}$ to the corresponding $\ol{a}_{\xi}$.
\end{claim}

\begin{cproof}
It suffices to verify that $F_X$ satisfies the universal property defining
$\xE(X)$, that is, for every \ppoag\ $\seq{G,e}$ with elements
$a_{\xi}\in G$ such that $0\leq a_{\xi}\leq e$, for $\xi\in X$, there
exists a (necessarily unique) positive homomorphism from $F_X$ to $G$ that
sends $\ol{e}$ to $e$ and each $\ol{a}_{\xi}$ to the
corresponding~$a_{\xi}$. This, in turn, amounts to verifying the following
statement:
 \begin{equation}\label{Eq:ImplOrd}
 m\ol{e}+\sum\famm{n_{\xi}\ol{a}_{\xi}}{\xi\in X}\geq0\ \Longrightarrow\ 
 me+\sum\famm{n_{\xi}a_{\xi}}{\xi\in X}\geq0,
 \end{equation}
for all $m\in\ZZ$ and all $\seqm{n_{\xi}}{\xi\in X}\in\ZZ^{(X)}$.
As $\seqm{n_{\xi}}{\xi\in X}$ has finite support, we may
assume without loss of generality that $X$ is finite. By
Claim~\ref{Cl:DefOrdol}, the premise of \eqref{Eq:ImplOrd} means that
$m+\sum\famm{n_{\xi}}{\xi\in Y}\geq0$ in $\ZZ$ for every $Y\in\Pow(X)$.
We shall conclude the proof by induction on~$|X|$. For $|X|=0$ it is
immediate. For $X=\set{\xi}$, $m\geq0$, and $m+n\geq0$, we compute
 \[
 me+na_{\xi}\geq me+(-m)a_{\xi}=m(e-a_{\xi})\geq 0.
 \]
Now the induction step. Pick $\eta\in X$, and set
$k=\max\set{0,-n_{\eta}}$. Hence
 \begin{equation}\label{Eq:Interpk}
 -n_{\eta}\leq k\leq m+\sum\famm{n_{\xi}}{\xi\in Y},
 \quad\text{for all }Y\subseteq X\setminus\set{\eta}.
 \end{equation}
Therefore, the element
 \[
 me+\sum\famm{n_{\xi}a_{\xi}}{\xi\in X}=(ke+n_{\eta}a_{\eta})+
 \left((m-k)e+\sum\famm{n_{\xi}a_{\xi}}{\xi\in X\setminus\set{\eta}}\right)
 \]
is, by the induction hypothesis, expressed as the sum of two elements of
$G^+$, thus it belongs to $G^+$.
\end{cproof}

It follows from Claim~\ref{Cl:GXfreeRel} that
 \begin{equation}\label{Eq:CharacIneq}
 m\be+\sum\famm{n_{\xi}\ba_{\xi}}{\xi\in X}\geq0\text{ if{f} }
 m+\sum\famm{n_{\xi}}{\xi\in Y}\geq0\text{ for all }Y\subseteq X,
 \end{equation}
for all $m\in\ZZ$ and all $\seqm{n_{\xi}}{\xi\in X}\in\ZZ^{(X)}$. Both
statements about unperforation and order-unit follow immediately.
\end{proof}

\begin{notation}\label{Not:FunctF}
We put $\xF=\xI\circ\xE$, the composition of the two functors $\xI$ and
$\xE$.
\end{notation}

By using Lemmas~\ref{L:DirLimFunc} and~\ref{L:PresDir}, we obtain that
$\xF$ is a direct limits preserving functor from the category of sets
(with maps) to the category of dimension groups (with positive
homomorphisms).

\begin{lemma}\label{L:jEembd}
The canonical map $j_{\xE(X)}\colon\xE(X)\to\xF(X)$ is an embedding, for
every set $X$.
\end{lemma}

\begin{proof}
We use the explicit description
of $\xE(X)$ given in the proof of Lemma~\ref{L:FunctE}. Denote by $B_X$
the additive group of all bounded maps from $\Pow(X)$ to $\ZZ$.
Observe, in particular, that $\ol{e}$ has index~$1$ in~$B_X$.
Hence, $\xE(X)\cong F_X$ embeds into the dimension group
$\seq{B_X,\ol{e}}$ with order-unit of index at most $1$. For any
interpolator~$\imath$ on~$B_X$, the structure $\seq{B_X,\imath}$ is a
member of $\cD$, in which $\xE(X)$ embeds.
\end{proof}

\emph{We shall always identify $\xE(X)$ with its image in $\xF(X)$}, so
that $j_{\xE(X)}$ becomes the inclusion map from $\xE(X)$ into $\xF(X)$.
Observe that despite what is suggested by the proof of
Lemma~\ref{L:jEembd}, the element $\be$ does not, as a rule, have index~$1$
in~$\xF(X)$, but~$2$. The reason for this discrepancy is that for
nonempty~$X$, the canonical map $g\colon\xF(X)\to B_X$ is not one-to-one,
even on the positive cone of $\xF(X)$.
Indeed, picking $\xi\in X$ and putting
$\bx=\ip(0,0,\ba_{\xi},\be-\ba_{\xi})$, we get
$\bx\in\xF(X)^+$. Furthermore, there exists a normalized positive
homomorphism $h\colon\seq{\xE(X),\be}\to\seq{\ZZ,2}$ such that
$h(\ba_{\xi})=1$ and there exists an interpolator $\imath$ on~$\ZZ$ such
that $\imath(0,0,1,1)=1$, so
$\fact{h}{\imath}(\bx)=\imath(0,0,h(\ba_{\xi}),h(\be-\ba_{\xi}))=
\imath(0,0,1,1)=1$, and so $\bx>0$.
However,
$2\bx\leq\ba_{\xi}+(\be-\ba_{\xi})=\be$, thus $2g(\bx)\leq g(\be)=\ol{e}$,
and so $g(\bx)=0$.

\section{Supports and subgroups in $\xF(X)$}\label{S:F(X)}

Throughout this section we shall fix a set $X$. For all $Y\subseteq X$,
we put $f_{Y,X}=\xI(e_{Y,X})$, the canonical embedding from $\xF(Y)$ into
$\xF(X)$. A \emph{support} of an element $\bx\in\xF(X)$ is a subset~$Y$
of~$X$ such that $\bx\in f_{Y,X}[\xF(Y)]$. As the functor $\xF$ preserves
direct limits, every element of $\xF(X)$ has a finite support.

Now put $s_{X,Y}=\xI(r_{X,Y})$, $\ol{r}_{X,Y}=e_{Y,X}\circ r_{X,Y}$,
and $\ol{s}_{X,Y}=\xI(\ol{r}_{X,Y})$.
Hence $\ol{r}_{X,Y}$ is an idempotent positive endomorphism of $\xE(X)$,
and it can be defined as in~\eqref{Eq:DefrXY}. Furthermore,
$s_{X,Y}\colon\xF(X)\to\xF(Y)$ while $\ol{s}_{X,Y}$ is
an idempotent positive endomorphism of $\xF(X)$.

\begin{lemma}\label{L:Idfs}
The following equations hold, for all $Y$, $Z\subseteq X$:
\begin{enumerate}
\item $f_{Y,X}\circ s_{X,Y}\circ f_{Y,X}=f_{Y,X}$.

\item $\ol{s}_{X,Y}\circ\ol{s}_{X,Z}=\ol{s}_{X,Y\cap Z}$.

\item $s_{X,Y}\circ f_{Z,X}=f_{Y\cap Z,Y}\circ s_{Z,Y\cap Z}$.
\end{enumerate}
\end{lemma}

\begin{proof}
Apply the functor $\xI$ to the following equations, whose verifications
are immediate (actually, it is easy to infer the first two equations from
the third one):
\begin{align*}
&e_{Y,X}\circ r_{X,Y}\circ e_{Y,X}=e_{Y,X};\\
&\ol{r}_{X,Y}\circ\ol{r}_{X,Z}=\ol{r}_{X,Y\cap Z};\\
&r_{X,Y}\circ e_{Z,X}=e_{Y\cap Z,Y}\circ r_{Z,Y\cap Z}.\tag*{\qed}
\end{align*}
\renewcommand{\qed}{}
\end{proof}

\begin{lemma}\label{L:EquivSupp}
Let $\bx\in\xF(X)$ and let $Y\subseteq X$. Then $Y$ is a support of $\bx$
if{f} \linebreak $\ol{s}_{X,Y}(\bx)=\bx$.
\end{lemma}

\begin{proof}
Suppose first that $\ol{s}_{X,Y}(\bx)=\bx$, and put $\by=s_{X,Y}(\bx)$.
Then $\bx=\ol{s}_{X,Y}(\bx)=f_{Y,X}(\by)$ belongs to $f_{Y,X}[\xF(Y)]$.
Conversely, suppose that $\bx=f_{Y,X}(\by)$, for some $\by\in\xF(Y)$.
Then, using Lemma~\ref{L:Idfs}(i), we obtain
 \begin{equation}
 \ol{s}_{X,Y}(\bx)=f_{Y,X}\circ s_{X,Y}\circ
 f_{Y,X}(\by)=f_{Y,X}(\by)=\bx.\tag*{\qed}
 \end{equation}
\renewcommand{\qed}{}
\end{proof}

\begin{corollary}\label{C:UniqSupp}
Every element of $\xF(X)$ has a least support, which is a finite subset
of~$X$.
\end{corollary}

\begin{proof}
Let $Y$ and $Z$ be supports of $\bx\in\xF(X)$. It follows
from Lemma~\ref{L:EquivSupp} and Lemma~\ref{L:Idfs}(ii) that
$\bx=\ol{s}_{X,Y}(\bx)=\ol{s}_{X,Z}(\bx)$, thus
$\bx=\ol{s}_{X,Y}\circ\ol{s}_{X,Z}(\bx)=\ol{s}_{X,Y\cap Z}(\bx)$, and so,
again by Lemma~\ref{L:EquivSupp}, $Y\cap Z$ is a support of $\bx$. As
$\bx$ has a finite support, the conclusion follows.
\end{proof}

We shall denote by $\supp(\bx)$ the least support of an element $\bx$ of
$\xF(X)$.

\begin{lemma}\label{L:RestrSupp}
Let $\bx\in\xF(X)$ and let $Y\subseteq X$. Then
$\supp(s_{X,Y}(\bx))\subseteq\supp(\bx)\cap Y$.
\end{lemma}

\begin{proof}
Put $Z=\supp(\bx)$. There is $\bz\in\xF(Z)$ such that $\bx=f_{Z,X}(\bz)$,
thus, using Lemma~\ref{L:Idfs}(iii), $s_{X,Y}(\bx)=s_{X,Y}\circ
f_{Z,X}(\bz)=f_{Y\cap Z,Y}\circ s_{Z,Y\cap Z}(\bz)$,
and so $s_{X,Y}(\bx)$ belongs to the image of $f_{Y\cap Z,Y}$.
\end{proof}

Now we shall define certain additive subgroups $G^X_Z$ of $\xF(X)$, for
$Z\in[X]^{\leq2}$. First, we put $G^X_{\es}=\ZZ\be$. Next, for any
$\xi\in X$, we denote by $G^X_{\set{\xi}}$ the subgroup of $\xF(X)$
generated by $\set{\ba_{\xi},\bb_{\xi}}$. Finally, for all distinct $\xi$,
$\eta\in X$, we put
 \[
 \bc_{\xi,\eta}=\ip(0,\ba_{\xi}+\ba_{\eta}-\be,\ba_{\xi},\ba_{\eta}),
 \]
and we denote by $G^X_{\set{\xi,\eta}}$ the subgroup of $\xF(X)$ generated
by $\set{\ba_{\xi},\ba_{\eta},\bb_{\xi},\bb_{\eta},\bc_{\xi,\eta}}$. As,
by Axiom (SYMM) (see Section~\ref{S:FunctI}),
$\bc_{\xi,\eta}=\bc_{\eta,\xi}$, this definition is correct. For
$\xi\in X$, we define a positive homomorphism
$\varphi_{\xi}\colon\ZZ^2\to G^X_{\set{\xi}}$, and for $\xi\neq\eta$ in
$X$, we define a positive homomorphism
$\psi_{\xi,\eta}\colon\ZZ^4\to G^X_{\set{\xi,\eta}}$, by the rules
 \begin{align}
 \varphi_{\xi}(x_0,x_1)&=x_0\ba_{\xi}+x_1\bb_{\xi},\label{Eq:defphi}\\
 \psi_{\xi,\eta}(x_0,x_1,x_2,x_3)&=
 x_0\bc_{\xi,\eta}+x_1(\ba_{\xi}-\bc_{\xi,\eta})
 +x_2(\ba_{\eta}-\bc_{\xi,\eta})
 +x_3(\bc_{\xi,\eta}+\be-\ba_{\xi}-\ba_{\eta}),
 \label{Eq:defpsi}
 \end{align}
for all $x_0$, $x_1$, $x_2$, $x_3\in\ZZ$.

\begin{lemma}\label{L:IntersGZ}\hfill
\begin{enumerate}
\item All maps $\varphi_{\xi}$, for $\xi\in X$, and $\psi_{\xi,\eta}$,
for $\xi\neq\eta$ in $X$, are isomorphisms.

\item $G^X_Y\cap G^X_Z=G^X_{Y\cap Z}$, for all $Y$, $Z\in[X]^{\leq2}$.
\end{enumerate}
\end{lemma}

\begin{proof}
By the definition of $\xE(X)$, there exists a unique positive
homomorphism $\tau_{\xi}\colon\xE(X)\to\ZZ^2$ that sends~$\be$
to $\seq{1,1}$, $\ba_{\xi}$ to $\seq{1,0}$, and $\ba_{\zeta}$ to
$\seq{0,0}$ for all $\zeta\in X\setminus\set{\xi}$. Fix any interpolator
$\imath$ on $\ZZ^2$ and set $\pi_{\xi}=\fact{(\tau_{\xi})}{\imath}$. Then
$\pi_{\xi}\circ\varphi_{\xi}$ fixes both vectors $\seq{1,0}$ and
$\seq{1,1}$, thus it is the identity. Therefore, $\varphi_{\xi}$ is an
embedding, and thus an isomorphism.

Now let $\xi\neq\eta$ in $X$. There exists a unique positive
homomorphism\linebreak
$\sigma_{\xi,\eta}\colon\xE(X)\to\ZZ^4$ such that
 \begin{align*}
 \sigma_{\xi,\eta}(\ba_{\xi})&=\seq{1,1,0,0},&
 \sigma_{\xi,\eta}(\ba_{\eta})&=\seq{1,0,1,0},\\
 \sigma_{\xi,\eta}(\be)&=\seq{1,1,1,1},&
 \sigma_{\xi,\eta}(\ba_{\zeta})&=\seq{0,0,0,0},\text{ for all }
 \zeta\in X\setminus\set{\xi,\eta}.
 \end{align*}
Let $\imath$ be any interpolator on $\ZZ^4$ and set
$\rho_{\xi,\eta}=\fact{(\sigma_{\xi,\eta})}{\imath}$. As
 \[
 \seq{0,0,0,0},\seq{1,0,0,-1}\leq\rho_{\xi,\eta}(\bc_{\xi,\eta})
 \leq\seq{1,1,0,0},\seq{1,0,1,0},
 \]
the only possibility is $\rho_{\xi,\eta}(\bc_{\xi,\eta})=\seq{1,0,0,0}$.
It follows that $\rho_{\xi,\eta}\circ\psi_{\xi,\eta}$ fixes each of the
vectors $\seq{1,0,0,0}$, $\seq{1,1,0,0}$, $\seq{1,0,1,0}$, and
$\seq{1,1,1,1}$, whence it is the identity. In particular,
$\psi_{\xi,\eta}$ is an embedding, but it is obviously surjective, thus
it is an isomorphism.

Now let $\xi\neq\eta$ in $X$, and let
$\bz\in G^X_{\set{\xi}}\cap G^X_{\set{\eta}}$. There are $x$, $y$, $x'$,
$y'\in\ZZ$ such that
$\bz=x\ba_{\xi}+y\bb_{\xi}=x'\ba_{\eta}+y'\bb_{\eta}$. Applying
$\rho_{\xi,\eta}$ yields $\seq{x,x,y,y}=\seq{x',y',x',y'}$, whence
$x=x'=y'=y$, and so $\bz=x\be\in G^X_{\es}$. Therefore, 
$G^X_{\set{\xi}}\cap G^X_{\set{\eta}}=G^X_{\es}$.

Finally, let $\xi$, $\eta$, $\zeta$ be distinct elements of $X$, and let
$\bz\in G^X_{\set{\xi,\eta}}\cap G^X_{\set{\xi,\zeta}}$. There are $x_i$,
$y_i\in\ZZ$, for $i<4$, such that
 \begin{equation}\label{Eq:2exprbz}
 \begin{aligned}
 \bz=x_0\bc_{\xi,\eta}+x_1(\ba_{\xi}-\bc_{\xi,\eta})
 +x_2(\ba_{\eta}-\bc_{\xi,\eta})
 +x_3(\bc_{\xi,\eta}+\be-\ba_{\xi}-\ba_{\eta})\\
 =y_0\bc_{\xi,\zeta}+y_1(\ba_{\xi}-\bc_{\xi,\zeta})
 +y_2(\ba_{\zeta}-\bc_{\xi,\zeta})
 +y_3(\bc_{\xi,\zeta}+\be-\ba_{\xi}-\ba_{\zeta}).
 \end{aligned}
 \end{equation}
{}From $\rho_{\xi,\eta}(\ba_{\zeta})=\seq{0,0,0,0}$ it follows that
$\rho_{\xi,\eta}(\bc_{\xi,\zeta})=\seq{0,0,0,0}$. Hence, applying
$\rho_{\xi,\eta}$ to \eqref{Eq:2exprbz} yields that
$\seq{x_0,x_1,x_2,x_3}=\seq{y_1,y_1,y_3,y_3}$, and thus $x_0=x_1$,
$x_2=x_3$, whence $\bz=x_0\ba_{\xi}+x_2\bb_{\xi}\in G^X_{\set{\xi}}$.
All other instances of~(ii) can be easily deduced from the two above.
\end{proof}

\section{Smoothening interpolators on $\xF(X)$}\label{S:Smooth}
In the present section we shall also fix a set $X$.

\begin{definition}\label{D:SmInterp}
An interpolator $\imath$ on $\xF(X)$ is \emph{smoothening of level $2$}, if
all subgroups $G^X_Z$ \pup{see Section~\ref{S:F(X)}}, for
$Z\in[X]^{\leq2}$, are closed under $\imath$.
\end{definition}

\begin{lemma}\label{L:SmInterp}
There exists a smoothening interpolator of level $2$ on $\xF(X)$.
\end{lemma}

\begin{proof}
For all $p=\seq{x,x',y,y'}\in\xF(X)^4$, we put
$\rng p=\set{x,x',y,y'}$. It follows from Lemma~\ref{L:IntersGZ}(ii) that
the set
 \[
 I(p)=\setm{Z\in[X]^{\leq2}}{\rng p\subseteq G^X_Z}
 \]
is closed under intersection, hence it has a greatest lower bound $Z_p$
in $\seq{\Pow(X),\subseteq}$,
which belongs to $I(p)$ in case $I(p)$ is nonempty (otherwise $Z_p=X$).
Put $H_p=\nobreak G^X_{Z_p}$, where we define $G^X_X=\xF(X)$.
So $H_p$ contains $\rng p$, and it follows from
Lemma~\ref{L:IntersGZ}(i) that~$H_p$ is a dimension group. Now we consider
the equivalence relation $\sim$ on $\xF(X)^4$ generated by all pairs
$\seq{x,x',y,y'}\sim\seq{x',x,y,y'}$ and
$\seq{x,x',y,y'}\sim\seq{x,x',y',y}$, for $x$, $x'$, $y$, $y'\in\xF(X)$,
and we pick a subset $C$ of $\xF(X)^4$ such that for each $p\in\xF(X)^4$
there exists a unique $\ol{p}\in C$ such that $p\sim\ol{p}$. For each
$p\in\xF(X)^4$, we put
 \[
 \ol{\imath}(p)=\begin{cases}
 \text{any }z\in H_p\text{ such that }x,x'\leq z\leq y,y',&
 \text{if }x,x'\leq y,y',\\
 0,&\text{otherwise},
 \end{cases}
 \]
and then we define $\imath(p)=\ol{\imath}(\ol{p})$, for all $p\in\xF(X)^4$.
Observe that if $\rng p\subseteq G^X_Z$, with $Z\in[X]^{\leq2}$, then
$Z_p\subseteq Z$, thus $H_p=G^X_{Z_p}\subseteq G^X_Z$, and thus
$\imath(p)\in H_p\subseteq G^X_Z$. Hence, all $G^X_Z$, for
$Z\in[X]^{\leq2}$, are closed under $\imath$. Therefore, $\imath$ is a
smoothening interpolator of level~$2$ on $\xF(X)$.
\end{proof}

\begin{lemma}\label{L:SmInterp2}
Let $\imath$ be a smoothening interpolator of level $2$ on $\xF(X)$. Then
for all $Z\in[X]^{\leq2}$ and all $\bx\in\xF(X)$ with support $Z$, the
element $\fact{(j_{\xE(X)})}{\imath}(\bx)$ belongs to~$G^X_Z$.
\end{lemma}

\begin{proof}
Put $g=j_{\xE(X)}$.
By the definition of a support, there exists $\bz\in\xF(Z)$ such that
$\bx=f_{Z,X}(\bz)$. Therefore, by using Lemma~\ref{L:Transiota},
 \[
 \fact{g}{\imath}(\bx)=\fact{g}{\imath}\circ f_{Z,X}(\bz)
 =\fact{g}{\imath}\circ\xI(e_{Z,X})(\bz)
 =\fact{(g\circ e_{Z,X})}{\imath}(\bz).
 \]
However, $\im(g\circ e_{Z,X})=\xE(Z)\subseteq G^X_Z$ and $G^X_Z$ is closed
under~$\imath$, hence, by Lemma~\ref{L:Hcliota}, the image of
$\fact{(g\circ e_{Z,X})}{\imath}$ is contained in $G^X_Z$. In particular,
using Lemma~\ref{L:Transiota}, we obtain that
$\fact{g}{\imath}(\bx)=\fact{(g\circ e_{Z,X})}{\imath}(\bz)$ belongs to
$G^X_Z$.
\end{proof}

\section{Proof of the main theorem}\label{S:ThmA}

Let $P(X,L,\mu)$ denote the following
statement:
\begin{quote}\em
$X$ is a set, $L$ is a lattice, and $\mu\colon\Dim L\onto\xF(X)^+$ is a
surjective \Vhom.
\end{quote}
We say that $\mu$ is \emph{zero-separating}, if $\mu^{-1}\set{0}=\set{0}$.

\begin{lemma}\label{L:RedLmod}
If $P(X,L,\mu)$ holds, then $P(X,L',\mu')$ holds for some
\emph{modular} lattice~$L'$ and some zero-separating $\mu'$.
\end{lemma}

\begin{proof}
It follows from Lemma~\ref{L:Vhomzs} that
$I=\setm{\bx\in\Dim L}{\mu(\bx)=0}$ is an o-ideal of $\Dim L$ and the map
$\ol{\mu}\colon(\Dim L)/I\onto\xF(X)^+$, $[\bx]_I\mapsto\mu(\bx)$ is a
\Vhom. However, it follows from Propositions~2.1
and~2.4 in \cite{WDim} that $(\Dim L)/I\cong\Dim(L/{\theta})$,
where~$\theta$ is the congruence of $L$ defined by $x\equiv_{\theta}y$
if{f} $\DD(x,y)\in I$, for all $x$,
$y\in L$. Hence, replacing $L$ by $L'=L/{\theta}$, it suffices to prove
that if $\mu$ separates zero, then $L$ is modular. If $\set{o,a,b,c,i}$
is a (possibly degenerate) pentagon of $L$, that is,
$o\leq c\leq a\leq i$, $a\wedge b=o$, and $b\vee c=i$, then
 \[
 \mu\DD(o,c)=\mu\DD(b,i)=\mu\DD(o,a)=\mu\DD(o,c)+\mu\DD(c,a),
 \]
thus, since $\xF(X)^+$ is cancellative, $\mu\DD(c,a)=0$. Therefore, since
$\mu$ separates zero, $\DD(c,a)=0$, and hence $a=c$. This proves
the modularity of $L$.
\end{proof}

Our main theorem is a consequence of the following more precise result.

\begin{theorem}\label{T:MainVNeg}
Let $X$ be a set, let $L$ a lattice, and let $\mu\colon\Dim L\onto\xF(X)^+$
be a \Vhom\ with image containing~$\be$. Then $|X|\leq\aleph_1$.
\end{theorem}

\begin{proof}
Suppose, to the contrary, that $|X|\geq\aleph_2$. It follows from
Lemma~\ref{L:RedLmod} that we may assume that $L$ is modular and
$\mu$ is zero-separating.

As $\mu$ is a monoid homomorphism and $\be\in\im\mu$, there are 
a natural number~$n$ and elements $u_i<v_i$ in~$L$, for $i<\nobreak n$,
such that $\be=\sum\famm{\mu\DD(u_i,v_i)}{i<n}$. For all $\xi\in X$, we
obtain, by applying refinement in $\xF(X)^+$ to the equation
 \[
 \ba_{\xi}+\bb_{\xi}=\sum\famm{\mu\DD(u_i,v_i)}{i<n},
 \]
decompositions of the form
 \begin{equation}\label{Eq:Decabxi}
 \ba_{\xi}=\sum\famm{\ba_{\xi,i}}{i<n},\quad
 \bb_{\xi}=\sum\famm{\bb_{\xi,i}}{i<n}
 \end{equation}
in $\xF(X)^+$ such that
 \begin{equation}\label{Eq:abxiiDuivi}
 \ba_{\xi,i}+\bb_{\xi,i}=\mu\DD(u_i,v_i),\text{ for all }i<n.
 \end{equation}
Since $L$ is modular and $\mu$ is a \Vhom, we are entitled to apply
Lemma~\ref{L:DecDD(a,b)} to the latter equation, and hence we obtain a
positive integer $\ell_{\xi,i}$ and a finite chain in~$L$ of the form
 \[
 u_i=x_{\xi,i}^0\leq x_{\xi,i}^1\leq\cdots\leq
 x_{\xi,i}^{2\ell_{\xi,i}}=v_i
 \]
such that
 \begin{align}
 \ba_{\xi,i}&=\sum\famm{\mu\DD(x_{\xi,i}^{2j},x_{\xi,i}^{2j+1})}
 {j<\ell_{\xi,i}},\label{Eq:axii1}\\
 \bb_{\xi,i}&=\sum
 \famm{\mu\DD(x_{\xi,i}^{2j+1},x_{\xi,i}^{2j+2})}
 {j<\ell_{\xi,i}}.\label{Eq:bxii1}
 \end{align}
Now we define
 \[
 \Phi(\xi)=\bigcup\famm{\supp\mu\DD(x_{\xi,i}^j,x_{\xi,i}^{j+1})}
 {i<n,\ j<2\ell_{\xi,i}},\quad\text{for all }\xi\in X.
 \]
By applying L\'az\'ar's Theorem (see Lemma~\ref{L:Laza}), we obtain a
subset $X_1$ of $X$ of cardinality $\aleph_2$ such that
 \begin{equation}\label{Eq:Laza}
 \eta\notin\Phi(\xi),\text{ for all distinct }\xi,\,\eta\in X_1.
 \end{equation}
By Lemma~\ref{L:SmInterp}, there exists a smoothening interpolator
$\imath$ of level~$2$ on $\xF(X_1)$. Now we put
 \begin{align}
 \pi&=\fact{(j_{\xE(X_1)})}{\imath}\circ s_{X,X_1},
 & \mu'&=\pi\circ\mu,\label{Eq:Defpi}\\
 \ba'_{\xi,i}&=\pi(\ba_{\xi,i}),&\bb'_{\xi,i}&=\pi(\bb_{\xi,i}),
 \qquad\text{for all }\xi\in X_1\text{ and all }i<n.
 \end{align}
For all $\xi\in X$, $i<n$, and $j<2\ell_{\xi,i}$, it follows from
Lemma~\ref{L:RestrSupp} that $\Phi(\xi)\cap X_1$ is a support of the
element $s_{X,X_1}\mu\DD(x_{\xi,i}^j,x_{\xi,i}^{j+1})$, hence, if
$\xi\in X_1$ and by using \eqref{Eq:Laza}, we obtain that $\set{\xi}$ is a
support of $s_{X,X_1}\mu\DD(x_{\xi,i}^j,x_{\xi,i}^{j+1})$. Therefore, by 
applying $\fact{(j_{\xE(X_1)})}{\imath}$ and using Lemma~\ref{L:SmInterp2},
we obtain
 \begin{equation}\label{Eq:mu'D(xxii)}
 \mu'\DD(x_{\xi,i}^j,x_{\xi,i}^{j+1})\in G^{X_1}_{\set{\xi}},
 \qquad\text{for all }\xi\in X_1,\ i<n,\text{ and }j<2\ell_{\xi,i}.
 \end{equation}
By applying $\pi$ to the equations \eqref{Eq:Decabxi}--\eqref{Eq:bxii1}
and observing that all elements of $\xE(X_1)$ are fixed under $\pi$,
we obtain the equations
 \begin{align}
 \ba_{\xi}&=\sum\famm{\ba'_{\xi,i}}{i<n}\text{ and }
 \bb_{\xi}=\sum\famm{\bb'_{\xi,i}}{i<n},
 \quad\text{for all }\xi\in X_1,\label{Eq:DecabxiX1}\\
 \ba'_{\xi,i}&=\sum\famm{\mu'\DD(x_{\xi,i}^{2j},x_{\xi,i}^{2j+1})}
 {j<\ell_{\xi,i}},\quad\text{for all }\xi\in X_1\text{ and all }i<n,
 \label{Eq:axii12}\\
 \bb'_{\xi,i}&=\sum
 \famm{\mu'\DD(x_{\xi,i}^{2j+1},x_{\xi,i}^{2j+2})}
 {j<\ell_{\xi,i}},\quad\text{for all }\xi\in X_1\text{ and all }i<n.
\label{Eq:bxii12}
 \end{align}
Fix $\xi\in X_1$ and $i<n$.
It follows from \eqref{Eq:mu'D(xxii)}, \eqref{Eq:axii12}, and
\eqref{Eq:bxii12} that both~$\ba'_{\xi,i}$ and~$\bb'_{\xi,i}$ belong to
$G^{X_1}_{\set{\xi}}$. However, it
follows from \eqref{Eq:DecabxiX1} that $0\leq\ba'_{\xi,i}\leq\ba_{\xi}$.
Since the isomorphism
$\varphi_{\xi}^{-1}\colon G^{X_1}_{\set{\xi}}\to\ZZ^2$ 
(see \eqref{Eq:defphi}) carries $\ba_{\xi}$ to $\seq{1,0}$, it follows that
 \begin{equation}\label{Eq:ba'0or1}
 \ba'_{\xi,i}\in\set{0,\ba_{\xi}}.
 \end{equation}
It follows from \eqref{Eq:ba'0or1},
\eqref{Eq:mu'D(xxii)}, and \eqref{Eq:axii12} that there exists
$j<\ell_{\xi,i}$ such that
 \begin{equation}\label{Eq:almall0axi}
 \mu'\DD(x_{\xi,i}^{2j'},x_{\xi,i}^{2j'+1})=0,\text{ for all }
 j'<\ell_{\xi,i}\text{ with }j'\neq j.
 \end{equation}
Similarly, $\bb'_{\xi,i}\in\set{0,\bb_{\xi}}$ and there exists
$k<\ell_{\xi,i}$ such that
 \begin{equation}\label{Eq:almall0bxi}
 \mu'\DD(x_{\xi,i}^{2k'+1},x_{\xi,i}^{2k'+2})=0,\text{ for all }
 k'<\ell_{\xi,i}\text{ with }k'\neq k.
 \end{equation}
We define an element $z_{\xi,i}\in L$ as follows:
 \[
 z_{\xi,i}=\begin{cases}
 x_{\xi,i}^{2j+1},&\text{if }j\leq k,\\
 x_{\xi,i}^{2k+2},&\text{if }j>k.
 \end{cases}
 \]
It follows easily from \eqref{Eq:axii12}, \eqref{Eq:bxii12},
\eqref{Eq:almall0axi}, and \eqref{Eq:almall0bxi} that the following
statements hold:
 \begin{align*}
 \ba'_{\xi,i}&=\mu'\DD(u_i,z_{\xi,i})\text{ and }
 \bb'_{\xi,i}=\mu'\DD(z_{\xi,i},v_i),\text{ if }j\leq k,\\
 \bb'_{\xi,i}&=\mu'\DD(u_i,z_{\xi,i})\text{ and }
 \ba'_{\xi,i}=\mu'\DD(z_{\xi,i},v_i),\text{ if }j>k.
 \end{align*}
Let $A(\xi,i)$ hold, if $\ba'_{\xi,i}=\mu'\DD(u_i,z_{\xi,i})$ and
$\bb'_{\xi,i}=\mu'\DD(z_{\xi,i},v_i)$, and let $B(\xi,i)$ hold, if
$\bb'_{\xi,i}=\mu'\DD(u_i,z_{\xi,i})$ and
$\ba'_{\xi,i}=\mu'\DD(z_{\xi,i},v_i)$.
What will matter for us is that the following property is satisfied:
 \begin{equation}\label{Eq:Charzxii}
 \text{Either }A(\xi,i)\text{ or }B(\xi,i)\text{ holds},
 \quad\text{for all }\xi\in X_1
 \text{ and all }i<n.
 \end{equation}
Now we denote by $U$ the powerset of $n=\set{0,1,\dots,n-1}$, and we put
 \[
 Y_u=\setm{\xi\in X_1}{(\forall i\in u)A(\xi,i)\text{ and }
 (\forall i\in n\setminus u)B(\xi,i)},\qquad\text{for all }u\in U.
 \]
\setcounter{claim}{0}
\begin{claim}\label{Cl:unionYu}
$X_1=\bigcup\famm{Y_u}{u\in U}$.
\end{claim}

\begin{cproof}
Let $\xi\in X_1$, and put $u=\setm{i<n}{A(\xi,i)}$. It follows from
\eqref{Eq:Charzxii} that $B(\xi,i)$ holds, for all $i\in n\setminus u$.
Therefore, $\xi\in Y_u$.
\end{cproof}

Now we put
$\bd_{\xi,\eta}=\sum\famm{\mu'\DD(z_{\xi,i},z_{\eta,i})}{i<n}$, for all
$\xi$, $\eta\in X_1$.

\begin{claim}\label{Cl:Ineqdxieta}
The following inequalities hold:
\begin{enumerate}
\item $\bd_{\xi,\zeta}\leq\bd_{\xi,\eta}+\bd_{\eta,\zeta}$, for all $\xi$,
$\eta$, $\zeta\in X_1$;

\item $\bd_{\xi,\eta}\leq\ba_{\xi}+\ba_{\eta},\bb_{\xi}+\bb_{\eta}$,
for all $u\in U$ and all $\xi$, $\eta\in Y_u$;

\item $\be\leq\ba_{\eta}+\bb_{\xi}+\bd_{\xi,\eta},
\ba_{\xi}+\bb_{\eta}+\bd_{\xi,\eta}$, for
all $u\in U$ and all $\xi$, $\eta\in Y_u$.
\end{enumerate}
\end{claim}

\begin{cproof}
Item (i) follows immediately from Lemma~\ref{L:DDDist}(iii).

Now let $u\in U$ and let $\xi$, $\eta\in Y_u$. Let $i<n$. If $i\in u$,
then, by using again Lemma~\ref{L:DDDist},
 \begin{align*}
 \mu'\DD(z_{\xi,i},z_{\eta,i})&\leq
 \mu'\DD(z_{\xi,i},u_i)+\mu'\DD(u_i,z_{\eta,i})=
 \ba'_{\xi,i}+\ba'_{\eta,i},\\
 \mu'\DD(u_i,v_i)&\leq\mu'\DD(u_i,z_{\eta,i})+
 \mu'\DD(z_{\eta,i},z_{\xi,i})+\mu'\DD(z_{\xi,i},v_i)\\
 &=\ba'_{\eta,i}+\bb'_{\xi,i}+\mu'\DD(z_{\xi,i},z_{\eta,i}),
 \end{align*}
while if $i\in n\setminus u$,
 \begin{align*}
 \mu'\DD(z_{\xi,i},z_{\eta,i})&\leq
 \mu'\DD(z_{\xi,i},v_i)+\mu'\DD(v_i,z_{\eta,i})=
 \ba'_{\xi,i}+\ba'_{\eta,i},\\
 \mu'\DD(u_i,v_i)&\leq\mu'\DD(u_i,z_{\xi,i})+
 \mu'\DD(z_{\xi,i},z_{\eta,i})+\mu'\DD(z_{\eta,i},v_i)\\
 &=\bb'_{\xi,i}+\ba'_{\eta,i}+\mu'\DD(z_{\xi,i},z_{\eta,i}),
 \end{align*}
so that in any case,
 \begin{align}
 \mu'\DD(z_{\xi,i},z_{\eta,i})&\leq\ba'_{\xi,i}+\ba'_{\eta,i},
 \label{Eq:dxietaileqaa}\\
 \mu'\DD(u_i,v_i)&\leq
 \ba'_{\eta,i}+\bb'_{\xi,i}+\mu'\DD(z_{\xi,i},z_{\eta,i}).
 \label{Eq:dxietaigeq}
 \end{align}
Symmetrically, we can obtain
 \begin{align}
 \mu'\DD(z_{\xi,i},z_{\eta,i})&\leq\bb'_{\xi,i}+\bb'_{\eta,i},
 \label{Eq:dxietaileqbb}\\
 \mu'\DD(u_i,v_i)&\leq
 \ba'_{\xi,i}+\bb'_{\eta,i}+\mu'\DD(z_{\xi,i},z_{\eta,i}),
 \label{Eq:dxietaigeq2}
 \end{align}
Adding together all inequalities
\eqref{Eq:dxietaileqaa}--\eqref{Eq:dxietaigeq2}, for $i<n$,
establishes both~(ii) and~(iii).
\end{cproof}

By Claim~\ref{Cl:unionYu}, there exists $u\in U$ such that
$|Y_u|=\aleph_2$. For the rest of the proof we fix such a subset~$u$.
We define $\Psi(\set{\xi,\eta})=\supp\bd_{\xi,\eta}$,
for all distinct $\xi$, $\eta\in Y_u$. Applying Kuratowski's Theorem (see
Lemma~\ref{L:Kura}) to the map~$\Psi$, we obtain distinct elements
$\alpha$, $\beta$, $\gamma\in Y_u$ such that
$\alpha\notin\Psi(\set{\beta,\gamma})$,
$\beta\notin\Psi(\set{\alpha,\gamma})$, and
$\gamma\notin\Psi(\set{\alpha,\beta})$.

Put $X_2=\set{\alpha,\beta,\gamma}$. It follows from Lemma~\ref{L:SmInterp}
that there exists a smoothening interpolator $\jmath$ of level~$2$ on
$\xF(X_2)$. Put $\pi'=\fact{(j_{\xE(X_2)})}{\jmath}\circ s_{X_1,X_2}$,
a positive homomorphism from $\xF(X_1)$ to $\xF(X_2)$. For
all distinct $\xi$, $\eta\in Y_u$, it follows from Lemma~\ref{L:RestrSupp}
that $\Psi(\set{\xi,\eta})\cap X_2$ is a support of the element
$s_{X_1,X_2}(\bd_{\xi,\eta})$. Hence, we obtain that the pair
$\set{\xi,\eta}$ is a support of
$s_{X_1,X_2}(\bd_{\xi,\eta})$, for all distinct $\xi$, $\eta\in X_2$.
Therefore, putting $\bd'_{\xi,\eta}=\pi'(\bd_{\xi,\eta})$, applying
$\fact{(j_{\xE(X_1)})}{\jmath}$, and using Lemma~\ref{L:SmInterp2}, we
obtain that
 \begin{equation}\label{Eq:projdxieta}
 \bd'_{\xi,\eta}\in G^{X_2}_{\set{\xi,\eta}},\quad\text{for all distinct }
 \xi,\,\eta\in X_2.
 \end{equation}
Applying the positive homomorphism $\pi'$ to the inequalities in
Claim~\ref{Cl:Ineqdxieta}, we obtain the following new inequalities, for
all distinct $\xi$, $\eta$, $\zeta\in X_2$:
\begin{enumerate}
\item $\bd'_{\xi,\zeta}\leq\bd'_{\xi,\eta}+\bd'_{\eta,\zeta}$

\item $\bd'_{\xi,\eta}\leq
\ba_{\xi}+\ba_{\eta},\bb_{\xi}+\bb_{\eta}$.

\item $\be\leq\ba_{\eta}+\bb_{\xi}+\bd'_{\xi,\eta},
\ba_{\xi}+\bb_{\eta}+\bd'_{\xi,\eta}$.
\end{enumerate}
By applying the isomorphism $\psi_{\xi,\eta}^{-1}$ (see
\eqref{Eq:defpsi}) to the inequalities (ii) and (iii) above, we obtain the
inequalities
 \begin{align*}
 \psi_{\xi,\eta}^{-1}(\bd'_{\xi,\eta})&\leq\seq{1,1,0,0}+\seq{1,0,1,0},\\
 \psi_{\xi,\eta}^{-1}(\bd'_{\xi,\eta})&\leq\seq{0,0,1,1}+\seq{0,1,0,1},\\
 \seq{1,1,1,1}&\leq\psi_{\xi,\eta}^{-1}(\bd'_{\xi,\eta})+
 \seq{1,0,1,0}+\seq{0,0,1,1},\\
 \seq{1,1,1,1}&\leq\psi_{\xi,\eta}^{-1}(\bd'_{\xi,\eta})+
 \seq{1,1,0,0}+\seq{0,1,0,1},
 \end{align*}
which leaves the only possibility
 \[
 \psi_{\xi,\eta}^{-1}(\bd'_{\xi,\eta})=\seq{0,1,1,0},
 \]
that is,
 \[
 \bd'_{\xi,\eta}=\ba_{\xi}+\ba_{\eta}-2\bc_{\xi,\eta}.
 \]
Therefore, applying the inequality (i) above
with $\seq{\xi,\eta,\zeta}=\seq{\alpha,\beta,\gamma}$, we obtain
 \begin{equation}\label{Eq:Finalabc}
 \bc_{\alpha,\beta}+\bc_{\beta,\gamma}\leq
 \ba_{\beta}+\bc_{\alpha,\gamma}
 \end{equation}
in $\xF(X_2)$. However, we shall now prove that \eqref{Eq:Finalabc} does
not hold. Indeed, the structure $\seq{\ZZ^2,\seq{2,1}}$ is a dimension
group with order-unit of index~$2$, thus it expands to some member
$\seq{\ZZ^2,\seq{2,1},\iota}$ of $\cD$ (where $\iota$ is an interpolator
on $\ZZ^2$) such that
 \[
 \iota(\seq{0,0},\seq{0,-1},\seq{1,0},\seq{1,0})=\seq{0,0}\text{ and }
 \iota(\seq{0,0},\seq{0,0},\seq{1,0},\seq{1,1})=\seq{1,0}.
 \]
Now there exists a unique normalized positive homomorphism
$h\colon\seq{\xE(X_2),\be}\to\seq{\ZZ^2,\seq{2,1}}$ such that
 \[
 h(\ba_{\alpha})=h(\ba_{\gamma})=\seq{1,0},
 \text{ and }h(\ba_{\beta})=\seq{1,1}.
 \]
By definition, $h(\be)=\seq{2,1}$, so we can compute
 \begin{align*}
 \fact{h}{\iota}(\bc_{\alpha,\gamma})&=
 \iota(\seq{0,0},h(\ba_{\alpha}+\ba_{\gamma}-\be),
 h(\ba_{\alpha}),h(\ba_{\gamma}))\\
  &=\iota(\seq{0,0},\seq{0,-1},\seq{1,0},\seq{1,0})\\
 &=\seq{0,0},\\
 \fact{h}{\iota}(\bc_{\alpha,\beta})&=
 \iota(\seq{0,0},h(\ba_{\alpha}+\ba_{\beta}-\be),
 h(\ba_{\alpha}),h(\ba_{\beta}))\\
  &=\iota(\seq{0,0},\seq{0,0},\seq{1,0},\seq{1,1})\\
 &=\seq{1,0},\\
 \fact{h}{\iota}(\bc_{\beta,\gamma})&=
 \iota(\seq{0,0},h(\ba_{\beta}+\ba_{\gamma}-\be),
 h(\ba_{\beta}),h(\ba_{\gamma}))\\
  &=\iota(\seq{0,0},\seq{0,0},\seq{1,1},\seq{1,0})\\
 &=\seq{1,0}.
 \end{align*}
Therefore, applying $\fact{h}{\iota}$ to the inequality
\eqref{Eq:Finalabc} yields the inequality $\seq{2,0}\leq\seq{1,1}$ (in
$\ZZ^2$!), a contradiction.
\end{proof}

\begin{corollary}\label{C:MainVNeg}
For any set $X$, the following conditions are equivalent:
\begin{enumerate}
\item There exists a lattice $L$ such that $\Dim L\cong\xF(X)^+$.

\item There exists a complemented modular lattice $L$ such that
$\Dim L\cong\xF(X)^+$.

\item There exists a von~Neumann regular ring $R$ such that
$V(R)\cong\xF(X)^+$.

\item There exists a locally matricial ring $R$ such that
$K_0(R)\cong\xF(X)$.

\item $|X|\leq\aleph_1$.
\end{enumerate}
\end{corollary}

\begin{proof}
(i)$\Rightarrow$(v) follows immediately from Theorem~\ref{T:MainVNeg}.

Now suppose that $|X|\leq\aleph_1$. Then $\xF(X)$ is a dimension group of
cardinality at most $\aleph_1$; moreover, it has an order-unit (namely,
$\be$). By \cite[Theorem~1.5]{GoHa86}, for any field $\FF$, there exists a
locally matricial algebra $R$ over $\FF$ such that $K_0(R)\cong\xF(X)$.
Hence (v) implies~(iv).

(iv)$\Rightarrow$(iii) is trivial, as $V(R)\cong K_0(R)^+$ for any locally
matricial ring (and, more generally, for any unit-regular ring) $R$.

Now assume (iii). Since $R$ is von~Neumann regular, it follows from
Lemma~\ref{L:V(R)DimL} that $V(R)\cong\Dim L$, where $L$ is the
(complemented modular) lattice of all principal right ideals of $M_2(R)$.
Hence $\Dim L\cong\xF(X)^+$, and so~(ii) holds.

Finally, (ii)$\Rightarrow$(i) is a tautology.
\end{proof}

We conclude the paper with a problem.

\begin{problem}
Is every conical refinement monoid of cardinality at most
$\aleph_1$ isomorphic to $\Dim L$, for some modular lattice $L$?
\end{problem}

Even for countable monoids the question above is open. It is formally
similar to the fundamental open problem raised by K.\,R. Goodearl
in his survey paper~\cite{Good95}, that asks which refinement monoids are
isomorphic to~$V(R)$ for some von~Neumann regular ring~$R$. A positive
answer to Goodearl's question would yield a positive answer to the problem
above, with $L$ sectionally complemented modular.


\begin{thebibliography}{99}
\bibitem{EHMR}
P. Erd\H os, A. Hajnal, A. M\'at\'e, and R. Rado,
``Combinatorial Set Theory: Partition Relations for Cardinals'',
Studies in Logic and the Foundations of Mathematics \textbf{106},
North-Holland, Amsterdam - New York - Oxford, 1984. 347~p.

\bibitem{Gpoag}
K.\,R. Goodearl,
``Partially Ordered Abelian Groups with Interpolation'',
Mathematical Surveys and Monographs, Vol. \textbf{20}. American
Mathematical Society, Providence,  R.I., 1986. xxii+336~p.

\bibitem{GvnRR}
K.\,R. Goodearl, ``Von Neumann Regular Rings'', Second edition, Krieger,
Malabar, Fl., 1991, xvi+412~p.

\bibitem{Good94}
K.\,R. Goodearl,
\emph{$K_0$ of regular rings with bounded index of nilpotence},
in: \emph{Abelian Group Theory and Related Topics}
(R. G\"obel, P. Hill, and W. Liebert, eds.),
Contemporary Mathematics \textbf{171} (1994), 173--199.

\bibitem{Good95}
K.\,R. Goodearl,
\emph{Von~Neumann regular rings and direct sum decomposition problems},
in: \emph{Abelian Groups and Modules, Padova 1994}
(A. Facchini and C. Menini, eds.), Dordrecht (1995) Kluwer, 249--255.

\bibitem{GoHa86}
K.\,R. Goodearl and D.\,E. Handelman,
\emph{Tensor products of dimension groups and $K_0$ of
unit-regular rings}, Canad. J. Math. \textbf{38}, no.~3 (1986), 633--658.

\bibitem{Gore}
K.\,R. Goodearl and F. Wehrung,
\emph{The complete dimension theory of partially ordered systems with
equivalence and orthogonality}, Mem. Amer. Math. Soc., to appear.

\bibitem{GLT2}
G. Gr\"atzer, ``General Lattice Theory. Second edition'', new appendices
by the author with B.\,A. Davey, R. Freese, B. Ganter, M. Greferath, P.
Jipsen, H.\,A. Priestley, H. Rose, E.\,T. Schmidt, S.\,E. Schmidt, F.
Wehrung, and R. Wille. Birkh\"auser Verlag, Basel, 1998. xx+663~p.

\bibitem{Kura51}
C. Kuratowski,
\emph{Sur une caract\'erisation des alephs},
Fund. Math. \textbf{38} (1951), 14--17.

\bibitem{Laza}
D. L\'az\'ar,
\emph{On a problem in the theory of aggregates}, Compositio
Math. \textbf{3} (1936), 304.

\bibitem{Malc}
A.\,I. Mal'cev,
``Algebraic Systems'',
Algebraic systems. (Algebraicheskie sistemy) (Russian)
Sovremennaja Algebra. Moskau: Verlag ``Nauka'', Hauptredaktion
f\"ur physikalisch-mathematische Literatur. 392 p., 1970.
English translation: Die Grundlehren der mathematischen Wissenschaften.
Band \textbf{192}. Berlin - Heidelberg - New York: Springer-Verlag; Berlin:
Akademie-Verlag, 1973. xii+317~p.

\bibitem{PTW}
M. Plo\v{s}\v{c}ica, J. T\r{u}ma, and F. Wehrung,
\emph{Congruence lattices of free lattices in nondistributive
varieties}, Colloq. Math. \textbf{76}, no.~2 (1998), 269--278.

\bibitem{TuWe}
J. T\r{u}ma and F. Wehrung,
\emph{Simultaneous representations of semilattices
by lattices with permutable congruences},
Internat. J. Algebra Comput. \textbf{11}, no.~2 (2001), 217--246.

\bibitem{CLPSurv}
J. T\r{u}ma and F. Wehrung,
\emph{A survey of recent results on congruence lattices of lattices},
Algebra Universalis \textbf{48}, no.~4 (2002), 439--471.

\bibitem{NonMeas}
F. Wehrung,
\emph{Non-measurability properties of interpolation vector spaces},
Israel J. Math. \textbf{103} (1998), 177--206.

\bibitem{WDim}
F. Wehrung,
\emph{The dimension monoid of a lattice}, Algebra Universalis
\textbf{40}, no.~3 (1998), 247--411.

\bibitem{UnifRef}
F. Wehrung,
\emph{A uniform refinement property for congruence lattices},
Proc. Amer. Math. Soc. \textbf{127}, no.~2 (1999), 363--370.


\end{thebibliography}
\end{document}